\documentclass[12pt, reqno]{amsart}
\usepackage{amsmath, amsthm, amscd, amsfonts, amssymb}
\usepackage[bookmarksnumbered]{hyperref}

\newtheorem{theorem}{Theorem}[section]

\theoremstyle{definition}

\newtheorem{example}[theorem]{Example}

\theoremstyle{remark}

\numberwithin{equation}{section}

\begin{document}

\title[Fuzzy stability of ring homomorphisms and ring derivations on ...]{Fuzzy stability of ring homomorphisms and ring derivations on fuzzy Banach algebras }
\author[M. Eshaghi Gordji and  N.Ghobadipour ]{M. Eshaghi Gordji and  N.Ghobadipour }
\address{Department of Mathematics, Semnan University, P. O. Box 35195-363, Semnan, Iran and Section of Mathematics and Informatics
} \email{madjid.eshaghi@gmail.com,    ghobadipour.n@gmail.com}



\begin{abstract} In
this paper, we establish  the Hyers--Ulam--Rassias stability of ring
homomorphisms and ring derivations on fuzzy Banach algebras.
\end{abstract}
\maketitle

%--------------------------------------------------------------------------------------%

\section{Introduction and preliminaries }
In 1984, Katsaras \cite{Kat} defined a fuzzy norm on a linear space
and at the same year Wu and Fang \cite{Co} also introduced a notion
of fuzzy normed space and gave the generalization of the Kolmogoroff
normalized theorem for  fuzzy topological linear space. In
\cite{Bi}, Biswas defined and studied fuzzy inner product spaces in
linear space. Since then some mathematicians have defined fuzzy
metrics and norms on a linear space from various points of view
\cite{Ba,Fe,Kri,Shi,Xi}. In 1994, Cheng and Mordeson introduced a
definition of fuzzy norm on a linear space in such a manner that the
corresponding induced fuzzy metric is of Kramosil and Michalek type
\cite{Kra}. In 2003, Bag and Samanta \cite{Ba} modified the
definition of Cheng and Mordeson \cite{Che} by removing a regular
condition. They also established a decomposition theorem of a fuzzy
norm into a family of crisp norms and investigated some properties
of fuzzy norms (see \cite{Bag2}). Following \cite{Bag1}, we give the
employing notion of a fuzzy
norm.\\
Let X be a real linear space. A function $N : X \times \Bbb R
\longrightarrow [0,1]$(the so-called fuzzy subset) is said to be a
fuzzy norm on X if for all $x, y \in X$ and all $a,b \in \Bbb
R$:\\
$(N_1)~~N(x,a)=0$ for $a\leq 0;$\\
$(N_2)~~x=0$ if and only if $N(x,a)=1$ for all $a>0;$\\
$(N_3)~~N(ax,b)=N(x,\frac{b}{|a|})$ if $a\neq0$;\\
$(N_4)~~N(x+y,a+b)\geq min \{N(x,a), N(y,b)\};$\\
$(N_5)~~N(x,.)$ is non-decreasing function on $\Bbb R$ and $\lim
_{a \to \infty} N(x,a)=1;$\\
$(N_6)~~$ For $x\neq 0,$ $N(x,.)$ is (upper semi) continuous on
$\Bbb R.$\\
The pair $(X,N)$ is called a fuzzy normed linear space. One may
regard $N(x,a)$ as the truth value of the statement "the norm of $x$
is less than or equal to the real number $a$".\\
\begin{example}
Let $(X,\|.\|)$ be a normed linear space. Then $$N(x,a)=\left\{%
\begin{array}{ll}
   \frac{a}{a+\|x\|}, & a>0~~~, x \in X, \\
    0,~~~ & a\leq0, x \in X \\
\end{array}%
\right.    $$ is a fuzzy norm on X.
\end{example}
Let $(X,N)$ be a fuzzy normed linear space. Let $\{x_n\}$ be a
sequence in X. Then $\{x_n\}$ is said to be convergent if there
exists $x \in X$ such that $\lim_{n \to \infty}N(x_n-x,a)=1$ for all
$a>0.$ In that case, $x$ is called the limit of the sequence
$\{x_n\}$ and we denote it by $N - lim_{n \to \infty} ~ x_n = x.$ A
sequence $\{x_n\}$ in X is called Cauchy if for each $\epsilon > 0$
and each $a_0$ there exists $n_0$ such that for all $n\geq n_0$ and
all $p > 0,$ we have $N(x_{n+p} - x_n, a) > 1$ - $\epsilon$. It is
known that every convergent sequence in  fuzzy normed space is
Cauchy. If each Cauchy sequence is convergent, then the fuzzy norm
is said to be complete and the fuzzy normed space is called a fuzzy
Banach space.\\ Let $X$ be an algebra and $(X,N)$ be complete fuzzy
normed space, the pair $(X,N)$ is said to be a fuzzy Banach algebra
if for every $x,y \in X,$ $a,b \in \Bbb R$  $$N(xy,ab)\geq min
\{N(x,a),N(y,b)\}.$$ As an example,
Let $(X,\|.\|)$ be a Banach algebra. Define, $$N(x,a)=\left\{%
\begin{array}{ll}
   0, & a\leq \|x\| \\
    1,~~~ & a>\|x\| \\
\end{array}%
\right.    $$ Then $(X,N)$ is a fuzzy Banach algebra.\\ The
stability problem of functional equations originated from a question
of Ulam \cite{Ul} in 1940, concerning the stability of group
homomorphisms. Let $(G_1,.)$ be a group and let $(G_2,*)$ be a
metric group with the metric $d(.,.).$ Given $\epsilon >0$, does
there exist a $\delta 0$, such that if a mapping
$h:G_1\longrightarrow G_2$ satisfies the inequality
$d(h(x.y),h(x)*h(y)) <\delta$ for all $x,y\in G_1$, then there
exists a homomorphism $H:G_1\longrightarrow G_2$ with
$d(h(x),H(x))<\epsilon$ for all $x\in G_1?$ In the other words,
under what condition does there exists a homomorphism near an
approximate homomorphism? The concept of stability for functional
equation arises when we replace the functional equation by an
inequality which acts as a perturbation of the equation. In 1941, D.
H. Hyers \cite{Hy1} gave the first affirmative  answer to the
question of Ulam for Banach spaces. Let $f:{E}\longrightarrow{E'}$
be a mapping between Banach spaces such that
$$\|f(x+y)-f(x)-f(y)\|\leq \delta $$
for all $x,y\in E,$ and for some $\delta>0.$ Then there exists a
unique additive mapping $T:{E}\longrightarrow{E'}$ such that
$$\|f(x)-T(x)\|\leq \delta$$
for all $x\in E.$ Moreover if $f(tx)$ is continuous in
$t\in\mathbb{R}$  for each fixed $x\in E,$ then $T$ is linear.
Finally in 1978, Th. M. Rassias \cite{TRa} proved the following
theorem.

\begin{theorem}\label{t1} Let $f:{E}\longrightarrow{E'}$ be a mapping from
 a normed vector space ${E}$
into a Banach space ${E'}$ subject to the inequality
$$\|f(x+y)-f(x)-f(y)\|\leq \epsilon (\|x\|^p+\|y\|^p) \eqno \hspace {0.5
 cm} (1.1)$$
for all $x,y\in E,$ where $\epsilon$ and p are constants with
$\epsilon>0$ and $p<1.$ Then there exists a unique additive
mapping $T:{E}\longrightarrow{E'}$ such that
$$\|f(x)-T(x)\|\leq \frac{2\epsilon}{2-2^p}\|x\|^p  \eqno \hspace {0.5
 cm}(1.2)$$ for all $x\in E.$
If $p<0$ then inequality $(1.1)$ holds for all $x,y\neq 0$, and
$(1.2)$ for $x\neq 0.$ Also, if the function $t\mapsto f(tx)$ from
$\Bbb R$ into $E'$ is continuous in real $t$ for each fixed $x\in
E,$ then $T$ is linear.
\end{theorem}
In 1991, Z. Gajda \cite{Gaj} answered the question for the case
$p>1$, which was raised by Rassias.  This new concept is known as
Hyers--Ulam--Rassias stability of functional equations.  The
stability of a generalization of the Bourgin's result on approximate
ring homomorphisms and derivations has proved by Roman Badora (see
\cite{Bad1,Bad2}). In the present paper, we investigate the
generalized Hyers--Ulam--Rassias stability of ring homomorphisms and
derivations in fuzzy normed spaces.
\section{Fuzzy stability of ring homomorphisms}
For ring homomorphisms Bourgin (see \cite{Bou}) proved the
following.
\begin{theorem}\label{t1}
Let $\epsilon$ and $\delta$ be nonnegative real numbers. Then
every mapping $f$ of a Banach algebra $X$ with an identity element
onto a Banach algebra $Y$ with an identity element satisfying
$$\|f(x+y)-f(x)-f(y)\|\leq \epsilon \hspace{9cm}$$ and $$\|f(x.y)-f(x)f(y)\|\leq
\delta,\hspace{9.5cm}$$for all $x,y \in X,$ is a ring homomorphism
of $X$ onto $Y,$ i.e., $$f(x+y)=f(x)+f(y)\hspace{9.5cm}$$ and
$$f(x.y)=f(x)f(y)\hspace{10cm}$$ for all $x,y \in X.$
\end{theorem}
Now, we investigate the fuzzy stability of ring homomorphisms.
\begin{theorem}\label{t2}
Let $A$ be a ring, $(B,N)$ be a fuzzy Banach algebra and $(C,N')$
be a fuzzy normed space. Let $\varphi: A \times A \to C$ be a
function such that for some $0<\alpha<2,$
$$N'(\varphi(2a,2b),t)\geq N'(\alpha \varphi(a,b),t) \eqno(2.1)$$
for all $a,b \in A$ and all $t>0.$ Suppose that $f: A \to B$ is a
function such that $$N(f(a+b)-f(a)-f(b),t)\geq
N'(\varphi(a,b),t)\eqno(2.2)$$ and $$N(f(a.b)-f(a)f(b),s)\geq
N'(\varphi(a,b),s)\eqno(2.3)$$ for all $a,b \in A$ and all
$t,s>0.$ Then there exists a unique ring homomorphism $h:A \to B$
such that
$$N(f(a)-h(a),t)\geq N'(\frac{2\varphi(a,a)}{2-\alpha},t),\eqno(2.4)$$ where
$a \in A$ and $t>0.$
\end{theorem}
\begin{proof}
Theorem $3.1$ of \cite{Mir} shows that there exists an additive
function $h: A \to B$ such that
$$N(f(a)-h(a),t)\geq N'(\frac{2\varphi(a,a)}{2-\alpha},t),$$  where
$a \in A$ and $t>0.$\\
Now we only need to show that $h$ is a multiplicative function.
Our inequality follows that
$$N(f(na)-h(na),t)\geq
N'(\frac{2\varphi(na,na)}{2-\alpha},t)\eqno(2.5)$$for all $a \in A$
and all $t>0.$ Thus $$N(n^{-1}f(na)-n^{-1}h(na),n^{-1}t)\geq
N'(\frac{2\varphi(na,na)}{2-\alpha},t)\eqno(2.6)$$for all $a \in A$
and all $t>0.$ By the additivity  of $h$ it is easy to see that then
$$N(n^{-1}f(na)-h(a),t)\geq
N'(\frac{2\varphi(na,na)}{2-\alpha},nt)\eqno(2.7)$$for all $a \in A$
and all $t>0.$ Sending $n$ to infinity in (2.7) and using $(N_2)$
and $(N_5),$ we see that $$h(a)=N-\lim_{n \to
\infty}n^{-1}f(na)\eqno(2.8)$$ for all $a \in A.$ Using inequality
$(2.3),$ we get
$$N(f((na).b)-f(na)f(b),s)\geq
N'(\varphi(na,b),s)\eqno(2.9)$$ for all $a,b \in A$ and all $s>0.$
Thus
$$N(n^{-1}[f((na).b)-f(na)f(b)],s)\geq
N'(\varphi(na,b),ns)\eqno(2.10)$$ for all $a,b \in A$ and all
$s>0.$ Sending $n$ to infinity in (2.10) and using $(N_2)$ and
$(N_5),$ we see that $$N-\lim_{n \to
\infty}n^{-1}[f((na)b)-f(na)f(b)]=0. \eqno(2.11)$$ Applying
$(2.8)$ and $(2.11)$ we have
\begin{align*}
h(a.b)&=N-\lim_{n \to \infty}n^{-1}f(n(a.b))=N-\lim_{n \to
\infty}n^{-1}f((na).b)\\
&=N-\lim_{n \to \infty}n^{-1}[f((na).b)-f((na).b)+f(na)f(b)]\\
&=h(a)f(b)
\end{align*}
for all $a,b \in A.$ The result of our calculation is the
following functional equation $$h(a.b)=h(a)f(b)\eqno(2.12)$$ for
all $a,b \in A.$ From this equation by the additivity of $h$ we
have
\begin{align*}
h(a)f(nb)=h(a.(nb))=h((na).b)=h(na)f(b)=nh(a)f(b)
\end{align*}
for all $a,b \in A.$ Therefore, $$h(a)n^{-1}f(nb)=h(a)f(b)$$ for
all $a,b \in A.$ Sending $n$ to infinity, by $(2.8),$ we see that
$$h(a)h(b)=h(a)f(b)\eqno(2.13)$$ for all $a,b \in A.$ Combining
this formula with equation $(2.12)$ we have that $h$ is a
multiplicative function.\\
To prove the uniqueness property of $h,$ assume that $h'$ is
another ring homomorphism satisfying $(2.4).$ Since both $h$ and
$h'$ are additive we deduce that
\begin{align*}
N(h(a)-h'(a),t)=N(h(na)-h'(na),nt)\geq
N'(\frac{2\varphi(na,na)}{2-\alpha},\frac{nt}{2})
\end{align*}
for all $a \in A$ and all $t >0.$ Letting $n$ to infinity we find
that $$N(h(a)-h'(a),t)=1$$ for all $a \in A$ and all $t >0.$ Hence
$$h(a)=h'(a)$$ for all $a \in A.$
\end{proof}
\section{Fuzzy stability of ring derivations}

Let $A$ be an algebra. A function $d: A \to A$ is called a ring
derivation if and only if it satisfies the following functional
equations $$d(a+b)=d(a)+d(b);$$  $$d(a.b)=a.d(b)+d(a).b$$ for all
$a,b \in A.$

In this section we prove the fuzzy stability of derivations in fuzzy
Banach algebras.
\begin{theorem}\label{t3}
Let $(A,N)$  be a fuzzy Banach algebra and $(B,N')$ be a fuzzy
normed space. Let $\varphi: A \times A \to B$ be a function such
that for some $0<\alpha<2,$
$$N'(\varphi(2a,2b),t)\geq N'(\alpha \varphi(a,b),t) \eqno(3.1)$$
for all $a,b \in A$ and all $t>0.$ Suppose that $f: A \to A$ is a
function such that $$N(f(a+b)-f(a)-f(b),t)\geq
N'(\varphi(a,b),t)\eqno(3.2)$$ and $$N(f(a.b)-af(b)-f(a)b,s)\geq
N'(\varphi(a,b),s)\eqno(3.3)$$ for all $a,b \in A$ and all
$t,s>0.$ Then there exists a unique ring derivation $d:A \to A$
such that
$$N(f(a)-d(a),t)\geq N'(\frac{2\varphi(a,a)}{2-\alpha},t),\eqno(3.4)$$ where
$a \in A$ and $t>0.$
\end{theorem}
\begin{proof}
Theorem $3.1$ of \cite{Mir} shows that there exists an additive
function $d: A \to A$ such that
$$N(f(a)-d(a),t)\geq N'(\frac{2\varphi(a,a)}{2-\alpha},t),\eqno(3.5)$$
for all $a \in A$ and $t>0.$ Now we only need to show that $d$
satisfies $$d(a.b)=a.d(b)+d(a).b$$ for all $a,b \in A.$ Our
inequality $(3.5)$ implies that
$$N(f(na)-d(na),t)\geq
N'(\frac{2\varphi(na,na)}{2-\alpha},t)\eqno(3.6)$$for all $a \in
A$ and all $t>0.$ Thus $$N(n^{-1}f(na)-n^{-1}d(na),n^{-1}t)\geq
N'(\frac{2\varphi(na,na)}{2-\alpha},t)\eqno(3.7)$$for all $a \in
A$ and all $t>0.$ By the additivity of $d$ it is easy to see that
$$N(n^{-1}f(na)-d(a),t)\geq
N'(\frac{2\varphi(na,na)}{2-\alpha},nt)\eqno(3.8)$$for all $a \in
A$ and all $t>0.$ Letting  $n$  to infinity in (3.8) and using
$(N_2)$ and $(N_5),$ we see that $$d(a)=N-\lim_{n \to
\infty}n^{-1}f(na)\eqno(3.9)$$ for all $a \in A.$ Using inequality
$(3.3),$ we get
$$N(f((na).b)-(na).f(b)f(na).b,s)\geq
N'(\varphi(na,b),s)\eqno(3.10)$$ for all $a,b \in A$ and all
$s>0.$ Condition $(3.3)$ implies that the function $g: A\times A
\to A$ defined by $$g(a.b)=f(a.b)-a.f(b)-f(a).b$$ is bounded for
all $a,b \in A.$ Hence, $$N-\lim_{n \to \infty}n^{-1}g(na,b)=0$$
for all $a,b \in A.$ Now, applying $(3.9)$ we get
$$d(a.b)=a.f(b)+d(a).b\eqno(3.11)$$ for all $a,b \in A.$ Indeed,
\begin{align*}
d(a.b)&=N-\lim_{n \to \infty}n^{-1}f(n(a.b)\\
&=N-\lim_{n \to \infty}n^{-1}f((na).b)\\
&=N-\lim_{n \to \infty}n^{-1}[na.f(b)+f(na).b+g(na,b)]\\
&=N-\lim_{n \to \infty}[a.f(b)+n^{-1}f(na).b+n^{-1}g(na,b)]\\
&=a.f(b)+d(a).b
\end{align*}
for all $a,b \in A.$ Let $a,b \in A$ and $n \in \Bbb N$ be fixed.
Then, using $(3.11)$ and the additivity of $d,$ we have
\begin{align*}
a.f(nb)+nd(a).b&=a.f(nb)+d(a).nb =d(a.nb)\\
&=d(na.b)=na.f(b)+d(na).b\\
&=na.f(b)+nd(a).b.
\end{align*}
Therefore, $$a.f(b)=a.n^{-1}f(nb)$$ for all $a,b \in A.$ Sending $n$
to infinity, by $(3.9),$ we see that
$$a.f(b)=a.d(b)\eqno(3.12)$$ for all $a,b \in A.$ Combining this
formula with equation $(3.11)$ we have
$$d(a.b)=a.d(b)+d(a).b$$ for all $a,b \in A.$
To prove the uniqueness property of $d,$ assume that $d'$ is
another ring derivation satisfying $(3.4).$ Since both $d$ and
$d'$ are additive we deduce that
\begin{align*}
N(d(a)-d'(a),t)=N(d(na)-d'(na),nt)\geq
N'(\frac{2\varphi(na,na)}{2-\alpha},\frac{nt}{2})
\end{align*}
for all $a \in A$ and all $t >0.$ Sending  $n$ to infinity we find
that $$N(d(a)-d'(a),t)=1$$ for all $a \in A$ and all $t >0.$ Hence
$$d(a)=d'(a)$$ for all $a \in A.$
\end{proof}

\end{document}